\documentclass[12pt]{article}
\usepackage{amsmath, amsfonts, amssymb, amsthm}
\usepackage[pdftex]{hyperref} 
\usepackage{url}
\usepackage[misc,geometry]{ifsym} 

\usepackage[ruled]{algorithm2e}

\SetAlFnt{\small}
\SetAlCapFnt{\small}
\SetAlCapNameFnt{\small}
\SetAlCapHSkip{0pt}
\IncMargin{-\parindent}

\numberwithin {equation}{section}

\newcommand{\R}{{\mathbb R}}
\newcommand{\E}{{\mathbb E}}
\newcommand{\Prob}{{\mathbb P}}
\def\reff#1{{\rm (\ref{#1})}}  
\def\myproof{{\bf Proof.\ }}
\def\e{\varepsilon}

\def\gdif{\mathrel{\mathop{-}\limits^{\smash{\raise-1ex\hbox{\,$*$\,}}}}}

\newcommand{\la}{\langle}
\newcommand{\ra}{\rangle}

\newtheorem{Lm}{Lemma}[section]

\newtheorem{Th}{Theorem}[section]

\theoremstyle{definition}

\title{Stochastic Intermediate Gradient Method for Convex Problems with Inexact Stochastic Oracle \footnote{The research presented in this paper was partially supported by Russian Foundation for Basic Research grants 14-01-00722-a, 15-31-20571-mol\_a\_ved.}}
\author{
	Pavel Dvurechensky 
  \and
	Alexander Gasnikov 
}

\begin{document}
\maketitle
\begin{abstract}
In this paper we introduce new methods for convex optimization problems with inexact stochastic oracle. First method is an extension of the intermediate gradient method proposed by Devolder, Glineur and Nesterov for problems with inexact oracle. 
Our new method can be applied to the problems with composite structure, stochastic inexact oracle and allows using non-Euclidean setup. We prove estimates for mean rate of convergence and probabilities of large deviations from this rate. 
Also we introduce two modifications of this method for strongly convex problems. For the first modification we prove mean rate of convergence estimates and for the second we prove estimates for large deviations from the mean rate of convergence.
All the rates give the complexity estimates for proposed methods which up to multiplicative constant coincide with lower complexity bound for the considered class of convex composite optimization problems with stochastic inexact oracle.
\end{abstract}
\section{Introduction}
In this paper we introduce new first-order methods for problems from rather wide class of convex composite optimization problems with stochastic inexact oracle. First order methods are widely developed since the earliest stage of development of optimization theory, see e.g. \cite{Evt_1982, Polyak_1987}.
The book \cite{NemYud_1983} started an activity in providing complexity bounds for proposed methods and for considered classes of problems (see also \cite{Nest_2004}). Later for convex problem with special structure ellipsoid methods (e.g. \cite{KhaTarErl_1988}) and interior-point methods \cite{NemNes_1994} were proposed. These methods have very fast convergence rate but have rather costly iterations requiring for solving $n$-dimensional problem about $n^3 - n^4$ \cite{Nest_2004} arithmetic operations. This makes them ineffective for large dimensions ($n > 10^5$). In last decade problems of large and huge dimension \cite{Nes_2012} became one of the main focus of the research in optimization methods due to large amount of application areas such as telecommunications, the Internet, traffic flows, machine learning, mechanical disign etc. Usually in this areas requirements for precision of the approximation of the optimal value are not very high. This allows to use first-order methods which converge slower but have nearly dimension independent rate of convergence and each their iteration requires about $n^2$ or less arithmetic operations. So it is important to develop new effective first-order methods. 

Let $E$ be a finite-dimensional real vector space and $E^*$ be its dual. We denote the value of linear function $g \in E^*$ at $x\in E$ by $\la g, x \ra$. Let $\|\cdot\|$ be some norm on $E$.
In this paper we consider {\it composite optimization} problem of the form
\begin{equation}
\min_{x \in Q} \{ \varphi(x) := f(x) + h(x)\},
\label{eq:PrStateInit} 
\end{equation}
where $Q \subset E$ is a closed convex set, $h(x)$ is a simple convex function, $f(x)$ is convex function with { \it stochastic inexact oracle}. This means that for every $x \in Q$ there are $f_{\delta,L}(x) \in \R$ and $g_{\delta,L}(x) \in E^*$ such that
\begin{equation}
0 \leq f(y) - f_{\delta,L}(x) - \la g_{\delta,L}(x) , y-x \ra \leq \frac{L}{2} \|x-y\|^2 + \delta, \quad \forall y \in Q,
\label{eq:dLOracle}
\end{equation}
and also that instead of $(f_{\delta,L}(x), g_{\delta,L}(x))$ (we will call this pair a $(\delta,L)$-oracle) we use their stochastic approximations $(F_{\delta,L}(x, \xi), G_{\delta,L}(x, \xi))$. 
The latter means that at all point $x \in Q$, we associate
with $x$ a random variable $\xi$ whose probability distribution is supported $\Xi \subset \R$ and such that
\begin{align}
& \E_{\xi} F_{\delta,L}(x, \xi) = f_{\delta,L}(x) \label{eq:expectF} \\
& \E_{\xi} G_{\delta,L}(x, \xi) = g_{\delta,L}(x) \label{eq:expectG} \\
& \E_{\xi} (\|G_{\delta,L}(x, \xi) - g_{\delta,L}(x)\|_*)^2 \leq \sigma^2. \label{eq:varG} 
\end{align}
Here $\|\cdot\|_*$ is the dual norm corresponding to $\|\cdot\|_E$: $\|g\|_* = \sup_{y \in E} \{ \la g,y \ra: \|y\|_E\leq 1 \}$.

Note that this class of problems is very wide and includes for example problems of stochastic optimization, smooth and non-smooth problems (see \cite{DGN_FOM_2011}), problems with error in gradient, such problems as LASSO \cite{Tib_1996}.

The work \cite{DGN_FOM_2011} considers the case of deterministic $(\delta,L)$-oracle. It is shown there that Dual Gradient Method for minimizing function $f(x)$ generates an approximate solution with the error $O\left(\frac{LR^2}{k} + \delta\right)$, where $R$ is the distance between the starting point and the solution and $k$ is the iteration counter. Also the authors show that Fast Gradient Method for minimizing function $f(x)$ generates an approximate solution with the error $O\left(\frac{LR^2}{k^2} + k\delta\right)$ and hence accumulates the error of the oracle. In \cite{DGN_IGM_2013} the same authors propose an Intermediate Gradient Method with the error $O\left(\frac{LR^2}{k^p} + k^{p-1}\delta\right)$, where $p\in [1,2]$ is chosen before the method starts. This method allows to choose the tradeoff between the rate of convergence and rate of error accumulation by choosing appropriate value of $p$. In the thesis \cite{Dev_PhD_2013} all the mentioned above methods are extended for non-Euclidean setup. Also in \cite{GhLan_SCO1_2012,GhLan_SCO2_2013} the authors construct the method for composite stochastic optimization which can be used both for smooth and for non-smooth functions, but they don't consider deterministic error of the oracle.  

In this paper we are in the framework of the mentioned above stochastic inexact oracle which means that on each step of the algorithm we get the information with both stochastic and deterministic error. This framework is more general than the one considered in \cite{GhLan_SCO1_2012,GhLan_SCO2_2013}. We generalize the Intermediate Gradient Method for the case of composite optimization problem and stochastic error. The result is Stochastic Intermediate Gradient Method (Algorithm \ref{alg:SIGM}) which can be used in non-Euclidean setup and has the mean rate of convergence of  $O\left(\frac{LR^2}{k^p} + \frac{\sigma R}{\sqrt{k}} + k^{p-1}\delta\right)$ (Theorem \ref{Th04}) which can be useful if the noise level $\sigma$ is rather small and the constant $L$ is large. Also with some so-called light-tail assumption about random variables $\xi$ we obtain the bound for large deviations from the mean rate (Theorem \ref{Th05}). Then we propose an accelerated method for strongly convex problems (Algorithm \ref{alg:SIGMA1}) and estimate its rate of convergence (Theorem \ref{Th06}).  Finally we introduce Algorithm \ref{alg:SIGMA2} which allows to control large deviations from the mean rate of convergence (Theorem \ref{Th07}).
It follows from the results of \cite{NemYud_1983}, \cite{Dev_PhD_2013} that the obtained mean rates of the above algorithms lead to complexity estimates which up to multiplicative constant coincide with lower complexity bounds for the considered class of convex composite optimization problems with stochastic inexact oracle.

\section{Notation and terminology}
We assume that we have chosen some norm $\|\cdot\|$ on $E$. By $\partial f(x)$ we denote subdifferential of the function $f(x)$ at point $x$. Also we need a {\it prox-function} $d(x)$ which is differential and strongly convex with the parameter 1 on $Q$ with respect to $\|\cdot\|$. Let $x_0$ be the minimizer of $d(x)$ on $Q$. By translating and scaling $d(x)$ if necessary, we can always ensure that
\begin{equation}
d(x_0) = 0, \quad d(x) \geq \frac12 \|x-x_0\|^2, \quad \forall x \in Q.
\label{eq:dDef}
\end{equation}
We define also the corresponding {\it Bregman distance}:
\begin{equation}
V (x, z) = d(x) - d(z) - \la \nabla d(z), x - z \ra.
\label{eq:BrDistDef}
\end{equation}
Due to the strong convexity of $d(x)$ with parameter 1, we have:
\begin{equation}
V (x, z) \geq \frac12 \|x-z\|^2, \quad \forall x,z \in Q.
\label{eq:VfromBelow}
\end{equation}

\section{Stochastic Intermediate Gradient Method}
Let $\{\alpha_i\}_{i \geq 0}$, $\{\beta_i\}_{i \geq 0}$, $\{B_i\}_{i \geq 0}$ be three sequences of coefficients satisfying 
\begin{align}
& \alpha_0 \in (0,1], \quad \beta_{i+1} \geq \beta_i > L , \quad \forall i \geq 0, \label{eq:alpbet1} \\
& 0\leq \alpha_i \leq B_i , \quad \forall i \geq 0, \label{eq:alpB} \\
& \alpha_k^2 \beta _k  \leq B_k \beta_{k-1} \leq \left( \sum_{i=0}^k{\alpha_i} \right) \beta_{k-1}, \quad \forall k \geq 1. \label{eq:alpbetB}
\end{align}
We define also $A_k = \sum_{i=0}^k{\alpha_i}$ and $\tau_{i}=\frac{\alpha_{i+1}}{B_{i+1}}$. Note that by definition $\alpha_{0} = A_0 = B_0$. The Stochastic Intermediate Gradient Method is described below as Algotithm \ref{alg:SIGM}. 

\begin{algorithm}[h!]
\SetAlgoNoLine
\KwIn{The sequences $\{\alpha_i\}_{i \geq 0}$, $\{\beta_i\}_{i \geq 0}$, $\{B_i\}_{i \geq 0}$, functions $d(x)$, $V(x,z)$.}
\KwOut{The point $y_{k}$.}

Compute $x_0 = \arg \min_{x\in Q} \{d(x)\}$.

Let $\xi_0$ be a realization of the random variable $\xi$.

Compute $G_{\delta,L}(x_0,\xi_0)$.

Compute 
\begin{equation}
y_0 = \arg \min_{x\in Q} \{ \beta_0 d(x) + \alpha_0 \la G_{\delta,L}(x_0,\xi_0), x-x_0 \ra + h(x)\}
\label{eq:algy0}
\end{equation}
$k$ = 0.

\Repeat{}{
Compute 
\begin{equation}
z_k = \arg \min_{x\in Q} \{ \beta_k d(x) + \sum_{i=0}^k{\alpha_i \la G_{\delta,L}(x_i,\xi_i), x-x_i \ra} + A_k h(x)\}
\label{eq:algzk}
\end{equation}
Let 
\begin{equation}
x_{k+1} = \tau_k z_k + (1- \tau_k) y_k
\label{eq:algxkp1}
\end{equation}
Let $\xi_{k+1}$ be a realization of the random variable $\xi$.

Compute $G_{\delta,L}(x_{k+1},\xi_{k+1})$.

Compute 
\begin{equation}
\hat {x}_{k+1} = \arg \min_{x\in Q} \{ \beta_k V(x,z_k) + \alpha_{k+1} \la G_{\delta,L}(x_{k+1},\xi_{k+1}), x-z_k \ra + \alpha_{k+1} h(x)\}.
\label{eq:alghxkp1}
\end{equation}
Let 
\begin{equation}
w_{k+1} = \tau_k \hat {x}_{k+1} + (1- \tau_k) y_k
\label{eq:algwkp1}
\end{equation}
Let 
\begin{equation}
y_{k+1} = \frac{A_{k+1}-B_{k+1}}{A_{k+1}} y_{k} + \frac{B_{k+1}}{A_{k+1}} w_{k+1}
\label{eq:algykp1}
\end{equation}
    }
\caption{Stochastic Intermediate Gradient Method}
\label{alg:SIGM}
\end{algorithm}

\subsection{General convergence rate}
Let us obtain the convergence rate of the proposed method in terms of the sequences $\{A_i\}$, $\{B_i\}$, and $\{\beta_i\}$. Denote by $\Psi _k (x) = \beta_k d(x) + \sum_{i=0}^k{\alpha_i \left[ F_{\delta,L}(x_i,\xi_i) +  \la G_{\delta,L}(x_i,\xi_i), x-x_i \ra + h(x) \right]} $, our {\it model} of the objective function, $\Psi_k^* = \min_{x \in Q} \Psi _k (x)$ its minimal
value on the feasible set and $\xi_{[k]} = (\xi_0, \dots, \xi_k)$ the history of the random process
after $k$ iterations.
Let us show that $\{y_k\}_{k\geq 0}$ and $\{\Psi_k(x)\}_{k\geq 0}$ define a sequence of estimate functions.

\begin{Lm}
For all $k\geq 0$ the following inequality holds
\begin{equation}
A_k \varphi (y_k) \leq \Psi_k^* + E_k
\label{eq:Lm1est1}
\end{equation}
with 
\begin{align}
&E_k = \sum_{i=0}^k{B_i \delta} + \sum_{i=0}^k{\frac{B_i}{\beta_i-L}\left( \|G_{\delta,L}(x_i,\xi_i)-g_{\delta,L}(x_i) \|_*\right)^2 } + \notag \\
& +\sum_{i=0}^k{\alpha_i (f_{\delta,L}(x_i) - F_{\delta,L}(x_i,\xi_i))} + \notag \\ 
& +\sum_{i=1}^k{(B_{i} - \alpha_{i})\frac{\alpha_{i}}{B_{i} }\la g_{\delta,L}(x_i)-G_{\delta,L}(x_i,\xi_i),z_{i-1}-y_{i-1} \ra }. \notag
\end{align}
\label{Lm01}
\end{Lm}

\myproof
Denote $f_i=f_{\delta,L}(x_i)$, $g_i=g_{\delta,L}(x_i)$, $F_i=F_{\delta,L}(x_i,\xi_i)$, $G_i=G_{\delta,L}(x_i,\xi_i)$.  
Note that for all $g \in E^*$, $x \in E$, $\zeta >0$:
\begin{equation}
\la g,x \ra + \frac{\zeta}{2} \|x\|^2 \geq -\frac{1}{\zeta} \|g\|_*^2.
\label{eq:usefulIneq}
\end{equation} 
Let us prove first that the statement is true for $k=0$.
\begin{equation}
\begin{array}{lcl}
\Psi_0^*& \stackrel{\reff{eq:algy0}}{=} &\beta_0 d(y_0) + \alpha_0 \left[F_0+ \la G_0, y_0-x_0 \ra + h(y_0) \right] \notag \\
& \stackrel{\reff{eq:dDef}}{\geq} &\frac{\beta_0}{2}\|y_0-x_0\|^2 + \alpha_0 \left[F_0+ \la G_0, y_0-x_0 \ra + h(y_0) \right] \notag \\
& \stackrel{\reff{eq:alpbet1}}{\geq} &\alpha_0 \left[F_0+ \la G_0, y_0-x_0 \ra + h(y_0) + \frac{\beta_0}{2}\|y_0-x_0\|^2\right] \notag \\
& = &\alpha_0 \left[f_0+ \la g_0, y_0-x_0 \ra + h(y_0) + \frac{L}{2}\|y_0-x_0\|^2\right] + \notag \\
&   &+ \alpha_0 \left[F_0-f_0+ \la G_0-g_0, y_0-x_0 \ra + \frac{\beta_0-L}{2}\|y_0-x_0\|^2\right] \notag \\
& \stackrel{\reff{eq:dLOracle},\reff{eq:usefulIneq}}{\geq} & \alpha_0 \left[f(y_0)+ h(y_0) - \delta \right] + \alpha_0 \left[F_0-f_0 \right] - \\
&& - \frac{\alpha_0}{\beta_0-L}\|G_0-g_0\|^2_*, \notag
\end{array}
\end{equation}
which is \reff{eq:Lm1est1} for $k=0$ since $\alpha_{0} = A_0 = B_0$.

Let us assume that \reff{eq:Lm1est1} is true for some $k \geq 0$ and prove that it is also true for $k +1$. 
Let $gh(z_k) \in \partial h(z_k)$. From the optimality condition in \reff{eq:algzk} we have:
\begin{equation}
\la \beta_k \nabla d(z_k) + \sum_{i=0}^k{\alpha_i G_i} + A_k gh(z_k), x-z_k \ra \geq 0 , \quad \forall x \in Q.
\label{eq:optCondzk}
\end{equation}

Using the inequality $\beta_{k+1} \geq \beta_k$ we get:
\begin{equation}
\begin{array}{lcl}
\Psi_{k+1}(x)& = &\beta_{k+1} d(x) + \sum_{i=0}^{k+1}{\alpha_i \left[F_i+ \la G_i, x-x_i \ra \right]} + A_{k+1} h(x)  \notag \\
& \stackrel{\reff{eq:BrDistDef}}{\geq} & \beta_{k} V(x,z_k) + \beta_{k}d(z_k) +\beta_{k} \la \nabla d(z_k), x - z_k \ra + \notag \\
& &+ \sum_{i=0}^{k+1}{\alpha_i \left[F_i+ \la G_i, x-x_i \ra \right]} + A_{k+1} h(x) \notag \\
& \stackrel{\reff{eq:optCondzk}}{\geq} & \beta_{k} V(x,z_k) + \beta_{k}d(z_k) + \sum_{i=0}^{k}{\alpha_i \left[F_i+ \la G_i, z_k-x_i \ra \right]} + \notag \\
&& +A_{k+1} h(x) + \la A_k gh(z_k), z_k-x \ra + \notag \\
&& + \alpha_{k+1} \left[F_{k+1}+ \la G_{k+1}, x-x_{k+1} \ra \right] \notag \\
& \geq & \beta_{k}d(z_k) + \sum_{i=0}^{k}{\alpha_i \left[F_i+ \la G_i, z_k-x_i \ra + h(z_k) \right]} + \notag \\
&& + \beta_{k} V(x,z_k) + \alpha_{k+1} \left[F_{k+1}+ \la G_{k+1}, x-x_{k+1} \ra + h(x) \right]\notag \\
&\stackrel{\reff{eq:algzk}}{=} & \Psi_{k}^*+\beta_{k} V(x,z_k) + \\ 
&& + \alpha_{k+1} \left[F_{k+1}+ \la G_{k+1}, x-x_{k+1} \ra + h(x) \right] 
\end{array}\label{eq:Lm1Proof1}
\end{equation}
In the last inequality we used that due to the convexity of $h(x)$ $A_{k+1} h(x) + \la A_k gh(z_k) , z_k -x \ra \geq A_k h(z_k) + \alpha_{k+1} h(x)$.

Also, since $A_k = (B_{k+1}-\alpha_{k+1}) + (A_{k+1}-B_{k+1})$  we have the following chain of inequalities.
\begin{equation}
\begin{array}{lcl}
 & &\Psi_{k}^*+ \alpha_{k+1} \left[F_{k+1}+ \la G_{k+1}, x-x_{k+1} \ra + h(x) \right] \notag \\
 &\stackrel{\reff{eq:Lm1est1}}{\geq} & A_k \varphi(y_k) - E_k + \alpha_{k+1} \left[F_{k+1}+ \la G_{k+1}, x-x_{k+1} \ra + h(x) \right] \notag \\
 & = & A_k h(y_k)+ (A_{k+1}-B_{k+1})f(y_k) + (B_{k+1}-\alpha_{k+1}) f(y_k) -  \notag \\
 & & -E_k + \alpha_{k+1} \left[F_{k+1}+ \la G_{k+1}, x-x_{k+1} \ra + h(x) \right] \notag \\
 &\stackrel{\reff{eq:dLOracle}}{\geq} & (A_{k+1}-B_{k+1})f(y_k) - E_k + \notag \\
&& + (B_{k+1}-\alpha_{k+1}) (f_{k+1} + \la g_{k+1} , y_k -x_{k+1} \ra ) +  \notag \\
&& + \alpha_{k+1} \left[F_{k+1}+ \la G_{k+1}, x-x_{k+1} \ra  \right] + \alpha_{k+1} h(x) + A_k h(y_k) \notag \\
 & = & (A_{k+1}-B_{k+1})f(y_k) - E_k + \alpha_{k+1} h(x) + A_k h(y_k) + \\
 & & + B_{k+1} F_{k+1} + (B_{k+1}-\alpha_{k+1}) (f_{k+1} - F_{k+1} ) + \\
 & & + (B_{k+1}-\alpha_{k+1}) \la g_{k+1} - G_{k+1}, y_k - x_{k+1} \ra + \\
 & & + \la G_{k+1}, (B_{k+1}-\alpha_{k+1}) (y_k - x_{k+1}) + \alpha_{k+1}(x-x_{k+1}) \ra \\
&\stackrel{\reff{eq:algxkp1}}{=} & B_{k+1} F_{k+1} + (B_{k+1}-\alpha_{k+1}) \left[ f_{k+1} - F_{k+1} + \right. \\
&& + \left. \la g_{k+1} - G_{k+1}, y_k - x_{k+1} \ra + h(y_k) \right] + \\
&& \alpha_{k+1} \la G_{k+1}, x-z_k \ra -E_k + \alpha_{k+1} h(x) + \\
&& + (A_{k+1}-B_{k+1})( f(y_k) + h(y_k)).
\end{array}
\end{equation}
Using \reff{eq:Lm1Proof1} we obtain:
\begin{equation}
\begin{array}{lcl}
\Psi_{k+1}^* & \geq & B_{k+1} F_{k+1} + \min_{x\in Q} \left\{ \beta_{k} V(x,z_k) + \right. \\
&& + \left. \alpha_{k+1} \la G_{k+1}, x-z_k \ra + \alpha_{k+1} h(x) \right\} - E_k\\
& & + (A_{k+1}-B_{k+1}) \varphi(y_k) +   \\
&& (B_{k+1}-\alpha_{k+1}) \left[ f_{k+1} - F_{k+1} +  \right. \\
&& \left. \la g_{k+1} - G_{k+1}, y_k - x_{k+1} \ra + h(y_k) \right] \\
&\stackrel{\reff{eq:alghxkp1}}{=}&  B_{k+1} F_{k+1} + \beta_{k} V(\hat{x}_{k+1},z_k) +  \\
&& +  \alpha_{k+1} \la G_{k+1}, \hat{x}_{k+1}-z_k \ra + \alpha_{k+1} h(\hat{x}_{k+1}) - E_k + \\
&& + (A_{k+1}-B_{k+1}) \varphi(y_k) + \\
&& + (B_{k+1}-\alpha_{k+1}) \left[ f_{k+1} - F_{k+1} + \right. \\
&& + \left. \la g_{k+1} - G_{k+1}, y_k - x_{k+1} \ra + h(y_k) \right] \\
&\stackrel{\reff{eq:VfromBelow}}{\geq}& B_{k+1} \bigl[ F_{k+1} + \tau_k \la G_{k+1}, \hat{x}_{k+1}-z_k \ra + \\
&& + \frac{\beta_k}{B_{k+1}} \|\hat{x}_{k+1}-z_k\|^2 \bigr] - \\
&& -E_k + (B_{k+1}-\alpha_{k+1}) \left[ f_{k+1} - F_{k+1} + \right. \\
&& + \left. \la g_{k+1} - G_{k+1}, y_k - x_{k+1} \ra  \right] + \\
&& B_{k+1} (\tau_k h(\hat{x}_{k+1}) + (1-\tau_k) h(y_k)) + \\
&& + (A_{k+1}-B_{k+1}) \varphi(y_k) \\
&\stackrel{\reff{eq:alpbetB},\reff{eq:algwkp1}}{\geq}& B_{k+1} \bigl[ F_{k+1} + \tau_k \la G_{k+1}, \hat{x}_{k+1}-z_k \ra + \\
&&  + \frac{\tau_k^2\beta_{k+1}}{2} \|\hat{x}_{k+1}-z_k\|^2 \bigr] - E_k + \\
&& + (B_{k+1}-\alpha_{k+1}) \left[ f_{k+1} - F_{k+1} + \right. \\
&& + \left. \la g_{k+1} - G_{k+1}, y_k - x_{k+1} \ra  \right] + B_{k+1}h(w_{k+1}) + \\
&& + (A_{k+1}-B_{k+1}) \varphi(y_k) \\
&\stackrel{\reff{eq:algxkp1},\reff{eq:algwkp1}}{\geq}& B_{k+1} \bigl[ F_{k+1} + \la G_{k+1}, w_{k+1}-x_{k+1} \ra + \\
&&  + \frac{\beta_{k+1}}{2} \| w_{k+1}-x_{k+1}\|^2 +  h(w_{k+1})\bigr] - E_k + \\
&& + (B_{k+1}-\alpha_{k+1}) \left[ f_{k+1} - F_{k+1} + \right. \\
&& + \left. \la g_{k+1} - G_{k+1}, y_k - x_{k+1} \ra  \right] + \\
&& + (A_{k+1}-B_{k+1}) \varphi(y_k)  \\
&=&  B_{k+1} \bigl[ f_{k+1} + \la g_{k+1}, w_{k+1}-x_{k+1} \ra + \\
&&  + \frac{L}{2} \| w_{k+1}-x_{k+1}\|^2 +  h(w_{k+1}) \bigr] - E_k + \\ 
&& + B_{k+1} \bigl[ F_{k+1} - f_{k+1} + \la G_{k+1} - g_{k+1} , w_{k+1}-x_{k+1} \ra + \\
&&  + \frac{\beta_{k+1}-L}{2} \| w_{k+1}-x_{k+1}\|^2 \bigr] + \\
&& + (B_{k+1}-\alpha_{k+1}) \left[ f_{k+1} - F_{k+1} + \right. \\
&& + \left. \la g_{k+1} - G_{k+1}, y_k - x_{k+1} \ra  \right] + \\
&& + (A_{k+1}-B_{k+1}) \varphi(y_k)  \\
&\stackrel{\reff{eq:dLOracle}}{\geq}& B_{k+1} (f (w_{k+1}) + h(w_{k+1}) - \delta) - E_k \\
&& + \alpha_{k+1} ( F_{k+1} - f_{k+1} ) +   \\
&& + (B_{k+1}-\alpha_{k+1})  \la g_{k+1} - G_{k+1}, y_k - x_{k+1} \ra + \\
&& + B_{k+1} \bigl[ \la G_{k+1} - g_{k+1} , w_{k+1}-x_{k+1} \ra + \\
&&  + \frac{\beta_{k+1}-L}{2} \| w_{k+1}-x_{k+1}\|^2 \bigr] + \\
&& + (A_{k+1}-B_{k+1}) \varphi(y_k)  \\
&\stackrel{\reff{eq:algykp1}, \reff{eq:usefulIneq}}{\geq}& A_{k+1} \varphi(y_{k+1}) -E_k - B_{k+1} \delta + \\
&&+ \alpha_{k+1} ( F_{k+1} - f_{k+1} ) + \\
&&+ (B_{k+1}-\alpha_{k+1})  \la g_{k+1} - G_{k+1}, y_k - x_{k+1} \ra - \\
&& - \frac{B_{k+1}}{\beta_{k+1}-L} \| g_{k+1} - G_{k+1}\|_*^2, 
\end{array} \notag
\end{equation}
which in view of \reff{eq:algxkp1} is \reff{eq:Lm1est1} for $k+1$.
\qed

\begin{Lm}
For all $k\geq 0$ the following inequality holds
\begin{equation}
\Psi_k(x) \leq  A_k \varphi (x) + \beta_k d(x) + \bar{E}_k(x), \quad \forall x \in Q
\label{eq:Lm2est1}
\end{equation}
with 
\begin{align}
&\bar{E}_k(x) = \sum_{i=0}^k \alpha_i \bigl[ F_{\delta,L}(x_i,\xi_i) - f_{\delta,L}(x_i) + \notag \\
& \la G_{\delta,L}(x_i,\xi_i) - g_{\delta,L}(x_i),x-x_{i} \ra \bigr]. \notag
\end{align}
\label{Lm02}
\end{Lm}

\myproof Using the notation introduced in Lemma \ref{Lm01} we have:
\begin{equation}
\begin{array}{lcl}
 \Psi_k(x) & = & \beta_k d(x) + \sum_{i=0}^k \alpha_i \bigl[  f_i + \la g_i,x-x_{i} \ra \bigr] + \\
&& + \sum_{i=0}^k \alpha_i \bigl[ F_i - f_i + \la G_i - g_i,x-x_{i} \ra \bigr] + A_k h(x) \\
&\stackrel{\reff{eq:dLOracle}}{\leq}& \beta_k d(x) + A_k \varphi(x) + \bar{E}_k(x).
\end{array} \notag
\end{equation}
\qed

Combining Lemma \ref{Lm01} and Lemma \ref{Lm02} we obtain the following result.
\begin{Th}
Assume that the function $f$ is endowed with a stochastic inexact oracle
with noise level $\sigma$, bias $\delta$ and constant $L$. Then the sequence $y_k$ generated by the Algorithm \ref{alg:SIGM}, when applied to the problem \reff{eq:PrStateInit}, satisfies
\begin{align}
& \varphi(y_k) - \varphi^* \leq \frac{1}{A_k}\biggl( \beta_k d(x^*) + \sum_{i=0}^k B_i \delta + \notag \\
& + \sum_{i=0}^k \frac{B_i}{\beta_i-L}  \|G_{\delta,L}(x_i,\xi_i)-g_{\delta,L}(x_i) \|_*^2 + \notag\\ 
& + \sum_{i=0}^k \alpha_i \la G_{\delta,L}(x_i,\xi_i) - g_{\delta,L}(x_i),x^*-x_{i} \ra  + \notag \\
& + \sum_{i=1}^k (B_{i} - \alpha_{i})\frac{\alpha_{i}}{B_{i} }\la G_{\delta,L}(x_i,\xi_i) - g_{\delta,L}(x_i),y_{i-1}-z_{i-1} \ra \biggr). \label{eq:Th1est}
\end{align}
\label{Th01}
\end{Th}

\myproof 
From the inequalities \reff{eq:Lm1est1}, \reff{eq:Lm2est1}, by the definition of $\Psi_k(x)$ and $\Psi_k^*$ we have:
$$
A_k \varphi(y_k) \leq \Psi_k^* + E_k \leq \Psi_k(x^*) + E_k \leq A_k \varphi^* + \beta_k d(x^*) + \bar{E}_k(x^*) + E_k,
$$  
which immediately gives the statement of the theorem.
\qed

\begin{Th}
Assume that the function $f$ is endowed with a stochastic inexact oracle
with noise level $\sigma$, bias $\delta$ and constant $L$. Then the sequence $y_k$ generated by the Algorithm \ref{alg:SIGM}, when applied to the problem \reff{eq:PrStateInit}, satisfies
\begin{align}
& \E_{\xi_0, \dots, \xi_k} \varphi(y_k) - \varphi^* \leq \frac{\beta_k d(x^*)}{A_k} + \frac{\sum_{i=0}^k B_i \delta}{A_k} +\notag \\
& + \frac{1}{A_k}
 \sum_{i=0}^k \frac{B_i}{\beta_i-L}  \sigma^2. \notag
\end{align}
\label{Th02}
\end{Th}
\myproof 
Since $\E_{\xi_i} \left[ G_i |\xi_{[i-1]} \right] = g_i$ and since $x_i$, $y_{i-1}$, and $z_{i-1}$ are deterministic functions of $(\xi_0, \dots, \xi_{i-1})$, we have $\E_{\xi_i} \left[ \la G_i-g_i , x^* - x_i \ra |\xi_{[i-1]} \right] = \E_{\xi_i} \left[ \la G_i-g_i , y_{i-1} - z_{i-1} \ra |\xi_{[i-1]} \right] = 0$. Therefore the expectation of fourth and fifth term in \reff{eq:Th1est} with respect to $ \xi_0, \dots, \xi_k $ is zero. Also by our assumption $\E_{\xi_i} \left[  \|G_i-g_i \|_*^2 |\xi_{[i-1]} \right] \leq \sigma^2$ and hence $\E_{ \xi_0, \dots, \xi_k} \left[\sum_{i=0}^k \frac{B_i}{\beta_i-L} \|G_i-g_i \|_*^2 \right] \leq \sum_{i=0}^k \frac{B_i}{\beta_i-L} \sigma^2 $. 
\qed

\subsection{General probability of large deviations}
In this section we obtain an upper bound on the probability of large deviation for the $\varphi(y_k) - \varphi^*$. To obtain our results we make the following additional assumptions
\begin{enumerate}
	\item $\xi_0, \dots, \xi_k$ are i.i.d random variables. \label{as1}
	\item $G_{\delta,L}(x,\xi)$ satisfies the light-tail condition
	$ \E_\xi \left[ \exp \left( \frac{\|G_{\delta,L}(x,\xi)-g_{\delta,L}(x) \|^2_*}{\sigma^2} \right) \right] \leq \exp(1)$. \label{as2}
	\item Set $Q$ is bounded with diameter $D = \max_{x,y\in Q}\|x~-~y\|$. \label{as3}
\end{enumerate} 
We will need the following lemmas.
\begin{Lm}[\cite{JudLanNemShap_2009}, \cite{Dev_PhD_2013}]
Let $\xi_0, \dots, \xi_k$ be a sequence of realizations of the i.i.d. random variables $X_0, \dots, X_k$ and let $\Delta_i = \Delta_i (\xi_{[i]})$ be a deterministic function of $\xi_{[i]}$ such that for all $i \geq 0$:
$$
\E \left[\exp \left( \frac{\Delta_i^2}{\sigma^2} \right)|\xi_{[i-1]} \right] \leq \exp(1)
$$
and $c_0, \dots, c_k$ is a sequence of positive coefficients. Then we have for any $k \geq 0$ and any $\Omega \geq 0$:
$$
\Prob \left( \sum_{i=0}^k  c_i \Delta_i^2 \geq (1+\Omega) \sum_{i=0}^k  c_i \sigma^2 \right) \leq \exp(-\Omega).
$$
\label{Lm03}
\end{Lm}
\begin{Lm}[\cite{LanNemShap_2012}, \cite{Dev_PhD_2013}]
Let $\xi_0, \dots, \xi_k$ be a sequence of realizations of the i.i.d. random variables $X_0, \dots, X_k$ and let $\Gamma_i$ and $\eta_i$  be a deterministic function of $\xi_{[i]}$ such that:
\begin{enumerate}
	\item $\E \left[\Gamma_i|\xi_{[i-1]} \right] = 0$.
	\item $|\Gamma_i| \leq c_i \eta_i $, where $c_i$ is positive deterministic constant.
	\item $\E \left[\exp \left( \frac{\eta_i^2}{\sigma^2} \right)|\xi_{[i-1]} \right] \leq \exp(1)$.
\end{enumerate}
Then for any $k \geq 0$ and any $\Omega \geq 0$:
$$
\Prob \left( \sum_{i=0}^k  \Gamma_i \geq \sqrt{3\Omega}\sigma \sqrt{\sum_{i=0}^k  c_i^2} \right) \leq \exp(-\Omega).
$$
\label{Lm04}
\end{Lm}

\begin{Th}
If the assumptions \ref{as1}, \ref{as2}, \ref{as3} are satisfied, then for all $k \geq 0$ and all $\Omega \geq 0$, the sequence generated by the SIGM satisfies:
\begin{align}
&\Prob \Bigg( \varphi(y_k) - \varphi^* \geq \frac{\beta_k d(x^*)}{A_k} + \frac{\sum_{i=0}^k B_i \delta}{A_k} +  \notag \\
& + \frac{1+\Omega}{A_k} \sum_{i=0}^k \frac{B_i}{\beta_i-L}  \sigma^2 + \frac{2D \sigma \sqrt{3\Omega}}{A_k} \sqrt{\sum_{i=0}^k \alpha_i^2} \Bigg) \leq 3\exp(-\Omega). \notag
\end{align}
\label{Th03}
\end{Th}

\myproof
From the Theorem \ref{Th01} we know that for the SIGM, the gap $\varphi(y_k) - \varphi^*$ can
be bounded from above by the sum of four quantities:
\begin{enumerate}
	\item deterministic $I_1(k) = \frac{\beta_k d(x^*)}{A_k} + \frac{\sum_{i=0}^k B_i \delta}{A_k}$,
	\item random  $I_2(k, \xi_{[k]}) = \frac{1}{A_k}\sum_{i=0}^k \frac{B_i}{\beta_i-L}  \|G_{\delta,L}(x_i,\xi_i)-g_{\delta,L}(x_i) \|_*^2$,
	\item random  $I_3(k, \xi_{[k]}) = \frac{1}{A_k} \sum_{i=1}^k (B_{i} - \alpha_{i})\frac{\alpha_{i}}{B_{i} }\la G_{\delta,L}(x_i,\xi_i) - g_{\delta,L}(x_i),y_{i-1}-z_{i-1} \ra$, 
	\item random  $I_4(k, \xi_{[k]}) = \frac{1}{A_k} \sum_{i=0}^k \alpha_i \la G_{\delta,L}(x_i,\xi_i) - g_{\delta,L}(x_i),x^*-x_{i} \ra$. 
\end{enumerate}
For $I_2(k, \xi_{[k]})$ using Lemma \ref{Lm03} with $\Delta_i = \|G_{\delta,L}(x_i,\xi_i)-g_{\delta,L}(x_i) \|_*$ and $c_i = \frac{B_i}{A_k(\beta_i-L)}$ we obtain: 
$$
\Prob \left( I_2(k, \xi_{[k]}) \geq \frac{1+\Omega}{A_k}
 \sum_{i=0}^k \frac{B_i}{\beta_i-L}  \sigma^2 \right) \leq \exp(-\Omega)
$$
for all $k \geq 0$ and $\Omega \geq 0$.

For $I_3(k, \xi_{[k]})$ using Lemma \ref{Lm04} with $\Gamma_i = (B_{i} - \alpha_{i})\frac{\alpha_{i}}{A_k B_{i} }\la G_{\delta,L}(x_i,\xi_i) - g_{\delta,L}(x_i),y_{i-1}-z_{i-1} \ra$, $\eta_i = \|G_{\delta,L}(x_i,\xi_i)-g_{\delta,L}(x_i) \|_*$ and $c_i = \frac{\alpha_i D}{A_k}$ we obtain: 
$$
\Prob \left( I_3(k, \xi_{[k]}) \geq \frac{D \sigma \sqrt{3\Omega}}{A_k} \sqrt{\sum_{i=1}^k \alpha_i^2}  \right) \leq \exp(-\Omega)
$$
for all $k \geq 0$ and $\Omega \geq 0$.

For $I_4(k, \xi_{[k]})$ using Lemma \ref{Lm04} with $\Gamma_i = \frac{\alpha_{i}}{A_k}\la G_{\delta,L}(x_i,\xi_i) - g_{\delta,L}(x_i),x^*-x_{i} \ra$, $\eta_i = \|G_{\delta,L}(x_i,\xi_i)-g_{\delta,L}(x_i) \|_*$ and $c_i = \frac{\alpha_i D}{A_k}$ we obtain: 
$$
\Prob \left( I_4(k, \xi_{[k]}) \geq \frac{D \sigma \sqrt{3\Omega}}{A_k} \sqrt{\sum_{i=0}^k \alpha_i^2}  \right) \leq \exp(-\Omega)
$$
for all $k \geq 0$ and $\Omega \geq 0$.
Combining these results we obtain the statement of the theorem.
\qed

\subsection{Choice of the coefficients}
In Theorem \ref{Th02} we have obtained mean rate of convergence for SIGM and in Theorem \ref{Th03} we have obtained bounds for probability of large deviations for the error of the algorithm. These results are formulated in terms of sequences $\{\alpha_i\}_{i \geq 0}$, $\{\beta_i\}_{i \geq 0}$, $\{B_i\}_{i \geq 0}$ satisfying \reff{eq:alpbet1}, \reff{eq:alpB}, \reff{eq:alpbetB}. Let us choose these sequences to obtain the rate of convergence of $\Theta\left(\frac{LR^2}{k^p} + \frac{\sigma R}{\sqrt{k}} + k^{p-1}\delta\right)$. Let $a \geq 1$ and $b \geq 0$ be some parameters. Let us assume that we know a number $R$ such that $ \sqrt{2d(x^*)} \leq R$. We choose 
\begin{align}
& \alpha_i = \frac{1}{a} \left(\frac{i+p}{p} \right)^{p-1}, \quad  \forall i \geq 0, \label{eq:chooseAlp}\\
& \beta_i = L + \frac{b\sigma}{R} (i+p+1)^{\frac{2p-1}{2}} , \quad  \forall i \geq 0, \label{eq:chooseBet}\\
& B_i = a \alpha_i^2 = \frac{1}{a} \left(\frac{i+p}{p} \right)^{2p-2}, \quad  \forall i \geq 0. \label{eq:chooseB}
\end{align}
Then inequalities \reff{eq:alpbet1} and \reff{eq:alpB} hold and we need to check that \reff{eq:alpbetB} also holds. 
Also we have 
\begin{equation}
A_k = \sum_{i=0}^k \alpha_i \geq \frac{1}{a} \int_0^k \left(\frac{x+p}{p} \right)^{p-1} dx + \alpha_0 \geq \frac{1}{a} \left(\frac{k+p}{p} \right)^{p}.
\label{eq:Akestim}
\end{equation}

Clearly for any $i \geq 0$: 
\begin{align}
& \alpha_k^2 = \frac{1}{a^2} \left(\frac{k+p}{p} \right)^{2p-2} \leq \frac{1}{a} \left(\frac{k+p}{p} \right)^{2p-2} \leq   \notag \\ 
& \frac{1}{a} \left(\frac{k+p}{p} \right)^{p} \leq A_k. \notag 
\end{align}
If we choose $a = 2^{\frac{2p-1}{2}}$ then
\begin{align}
& \frac{1}{a^2} \left(\frac{k+p}{p} \right)^{2p-2} (k+p+1)^{\frac{2p-1}{2}} \leq  \notag  \\
& \leq \frac{1}{a} \left(\frac{k+p}{p} \right)^{2p-2} (k+p)^{\frac{2p-1}{2}} \leq \frac{1}{a} \left(\frac{k+p}{p} \right)^{p} (k+p)^{\frac{2p-1}{2}} . \notag 
\end{align}
Last two sequences of inequalities prove that \reff{eq:alpbetB} holds.

Using \reff{eq:Akestim} we have
\begin{align}
& \frac{\beta_k d(x^*)}{A_k} \leq \frac{\beta_k R^2}{2 A_k} \leq \left(  L + \frac{b\sigma}{R} (k+p+1)^{\frac{2p-1}{2}} \right) R^2 2^{\frac{2p-3}{2}} \left(\frac{p}{k+p} \right)^{p}.\label{eq:1termest}
\end{align}

Also using \reff{eq:Akestim} and the fact that $p \in [1,2]$ we have the following chain of inequalities
\begin{align}
& \frac{\delta}{A_k} \sum_{i=0}^k B_i = \frac{a \delta}{A_k}\sum_{i=0}^k \alpha_i^2 \leq \notag \\ 
&\leq \frac{a \delta}{A_k} \left( \int_0^k{ \left(\frac{x+p}{p} \right)^{2p-2} dx } + \left(\frac{k+p}{p} \right)^{2p-2} \right) \leq \notag \\ 
& \leq \frac{a \delta}{A_k} \left( \left(\frac{k+p}{p} \right)^{2p-1} + \left(\frac{k+p}{p} \right)^{2p-2} \right) \leq \notag \\
& \leq 2^{2p-1} \delta \left(\frac{p}{k+p} \right)^{p} \left( \left(\frac{k+p}{p} \right)^{2p-1} + \left(\frac{k+p}{p} \right)^{2p-2} \right) \leq \notag \\
& \leq 2^{2p-1} \left( \left(\frac{k+p}{p} \right)^{p-1} + 1 \right) \delta  .  \label{eq:2termest}
\end{align}
Using \reff{eq:Akestim} we have the following inequalities
\begin{align}
& \frac{\sigma^2}{A_k} \sum_{i=0}^k \frac{B_i}{\beta_i-L}  \leq \frac{\sigma R}{b p^{2p-2}} \left(\frac{p}{k+p} \right)^{p} \sum_{i=0}^k \frac{(i+p)^{2p-2}}{(i+p+1)^{\frac{2p-1}{2}}} \leq \notag \\
& \leq \frac{\sigma R p^{2 - p}}{b (k+p)^{p}} \sum_{i=0}^k {(i+p+1)^{p-\frac{3}{2}}} \leq \notag \\
& \leq \frac{\sigma R p^{2 - p}}{b (k+p)^{p}} \int_1^{k+1}{(x+p+1)^{p-\frac{3}{2}}dx} \leq \notag \\
& \leq \frac{\sigma R p^{2 - p}}{b (p-\frac12) }   \frac{(k+p+2)^{p-\frac12}}{(k+p)^{p}}. \label{eq:3termest}  
\end{align}
Combining estimates \reff{eq:1termest}, \reff{eq:2termest}, \reff{eq:3termest} we get for the estimation from Theorem \ref{Th02}:
\begin{align}
& \E_{X_0, \dots, X_k} \varphi(y_k) - \varphi^* \leq \notag \\
& \left(  L + \frac{b\sigma}{R} (k+p+1)^{\frac{2p-1}{2}} \right) R^2 2^{\frac{2p-3}{2}} \left(\frac{p}{k+p} \right)^{p} + \notag \\
& + \frac{\sigma R p^{2 - p}}{b (p-\frac12) }   \frac{(k+p+2)^{p-\frac12}}{(k+p)^{p}} + \notag \\ 
& + 2^{2p-1} \left( \left(\frac{k+p}{p} \right)^{p-1} + 1 \right) \delta \leq \notag \\
& \leq \frac{LR^2 p^p 2^{\frac{2p-3}{2}} }{(k+p)^{p}} + \notag \\  
& +\frac{\sigma R (k+p+2)^{p-\frac12}}{(k+p)^{p}} \left(b 2^{p-\frac32} p^p + \frac{2p^{1-p}}{b} \right) +  \notag \\
& + 2^{2p-1} \left( \left(\frac{k+p}{p} \right)^{p-1} + 1 \right) \delta. \notag
\end{align}
Choosing optimal $b=2^{\frac{5-2p}{4}}p^{\frac{1-2p}{2}}$ we get the following theorem.
\begin{Th}
If the sequences $\{\alpha_i\}_{i \geq 0}$, $\{\beta_i\}_{i \geq 0}$, $\{B_i\}_{i \geq 0}$ are chosen from relations \reff{eq:chooseAlp}, \reff{eq:chooseBet}, \reff{eq:chooseB} with $a = 2^{\frac{2p-1}{2}}$ and $b=2^{\frac{5-2p}{4}}p^{\frac{1-2p}{2}}$ then the
sequence $y_k$ generated by the SIGM satisfies:
\begin{align}
& \E_{X_0, \dots, X_k} \varphi(y_k) - \varphi^* \leq \notag \\
& \leq \frac{LR^2 p^p 2^{\frac{2p-3}{2}} }{(k+p)^{p}} + \frac{\sigma R 2^{\frac{3+2p}{4}  }\sqrt{p}(k+p+2)^{p-\frac12}}{(k+p)^{p}} +  \notag \\
& + 2^{2p-1} \left( \left(\frac{k+p}{p} \right)^{p-1} + 1 \right) \delta \leq \frac{C_1 LR^2  }{k^{p}} + \frac{C_2 \sigma R }{\sqrt{k}} +  C_3 k^{p-1}\delta =  \notag \\
& = \Theta\left(\frac{LR^2}{k^p} + \frac{\sigma R}{\sqrt{k}} + k^{p-1}\delta\right) ,\notag
\end{align}
where $C_1=4 \sqrt{2}$, $C_2 = 16 \sqrt{2}$, $C_3 = 48$.
\label{Th04}
\end{Th}

Similarly to what we have done to prove \reff{eq:2termest} we can get the following inequality:
$$
\frac{1}{A_k^2} \sum_{i=0}^k \alpha_i^2 \leq \frac{2p}{k+p}.
$$

This with \reff{eq:1termest}, \reff{eq:2termest}, \reff{eq:3termest} gives us the following corollary of Theorem \ref{Th03}.
\begin{Th}
If the sequences $\{\alpha_i\}_{i \geq 0}$, $\{\beta_i\}_{i \geq 0}$, $\{B_i\}_{i \geq 0}$ are chosen from relations \reff{eq:chooseAlp}, \reff{eq:chooseBet}, \reff{eq:chooseB} with $a = 2^{\frac{2p-1}{2}}$ and $b=2^{\frac{5-2p}{4}}p^{\frac{1-2p}{2}}$ then the
sequence $y_k$ generated by the SIGM satisfies:
\begin{align}
& \Prob \Biggl( \varphi(y_k) - \varphi^* > \notag \\
& > \frac{C_1 LR^2 }{k^{p}} + \frac{C_2( 1+ \Omega )\sigma R }{\sqrt{k}} +  \notag \\
& + C_3 k^{p-1}  \delta + \frac{C_4 D \sigma \sqrt{\Omega} }{\sqrt{k}} \Biggr) \leq  \notag \\
& \leq \Prob \Biggl( \varphi(y_k) - \varphi^* > \notag \\
& > \frac{LR^2 p^p 2^{\frac{2p-3}{2}} }{(k+p)^{p}} + \frac{( 1+ \Omega )\sigma R 2^{\frac{3+2p}{4}  }\sqrt{p}(k+p+2)^{p-\frac12}}{(k+p)^{p}} +  \notag \\
& + 2^{2p-1} \left( \left(\frac{k+p}{p} \right)^{p-1} + 1 \right) \delta + \frac{2D \sigma \sqrt{6 \Omega p }}{\sqrt{k+p}} \Biggr) \notag \leq \\ 
& \leq 3 \exp(- \Omega)  , \notag
\end{align}
where $C_1=4 \sqrt{2}$, $C_2 = 16 \sqrt{2}$, $C_3 = 48$, $C_4 = 4 \sqrt{3}$.
\label{Th05}
\end{Th}

\section{Stochastic Intermediate Gradient Method Accelerated}
In this section we consider two modifications of the SIGM method for strongly convex problems. For the first modification we obtain mean rate of convergence and for the second we bound the probability of large deviations from this rate. Both modifications are based on the restart technique which was previously used in \cite{GhLan_SCO2_2013} and \cite{JudNest_2014}.

Throughout this section we assume that $E$ is Euclidean space with scalar product $\la \cdot , \cdot \ra$ and norm $\|x\| = \sqrt{\la x , H x \ra}$, where $H$ is symmetric positive definite matrix.
Also we assume that the function $\varphi(x)$ in \reff{eq:PrStateInit} is strongly convex: 
$$
\frac{\mu}{2}\|x-y\|^2 \leq \varphi(y) - \varphi(x) - \la g(x) , y-x \ra, \quad \forall x,y \in Q, \quad g(x) \in \partial \varphi(x).
$$
As a corollary we have 
\begin{equation}
\varphi(x)-\varphi(x^*) \geq \frac{\mu}{2} \|x-x^*\|^2, \quad \forall x \in Q,
\label{eq:PhiStrConv}
\end{equation}
where $x^*$ is the solution of the problem \reff{eq:PrStateInit}.

Let us assume without loss of generality that the function $d(x)$ satisfies conditions $0=\arg \min_{x\in Q} d(x)$ and $d(0)=0$. 

\subsection{Method with Mean Rate of Convergence}
In this subsection we make the following assumption on the prox-function $d(x)$. We assume that if $x_0$ is random vector such that $\E_{x_0} \|x-x_0\|^2 \leq R_0^2$ for some fixed point $x$ and number $R_0$ then 
\begin{equation}
\E_{x_0} d \left( \frac{x-x_0}{R_0} \right) \leq \frac{V^2}{2}
\label{eq:dRestr}
\end{equation}
for some $V >0$. This assumption is satisfied for example for prox-functions with quadratic growth with constant $V^2$ which means that $d(x) \leq \frac{V^2}{2}\|x\|^2$ for all $x \in E$. Several examples of such prox-functions can be found in \cite{JudNest_2014}. Using this assumption we can obtain the following corollary of the Theorem \ref{Th04}.

\begin{Lm}
Assume that we start the Algorithm \ref{alg:SIGM} from random point $x_0$ such that $\E_{x_0} \|x^*-x_0\|^2 \leq R_0^2$ and hence \reff{eq:dRestr} holds with $x=x^*$. We use the function $d \left( \frac{x-x_0}{R_0} \right)$ as the prox-function in the algorithm. Also assume that on $k$-th iteration of the Algorithm \ref{alg:SIGM} we ask oracle $m$ times, getting answers $G_{\delta,L}(x_{k+1},\xi_{k+1}^i), \quad i=1,\dots,m$ and use $\tilde{G}_{\delta,L}(x_{k+1}) = \frac{1}{m}\sum_{i=1}^m{G_{\delta,L}(x_{k+1},\xi_{k+1}^i)}$ in \reff{eq:alghxkp1} instead of $G_{\delta,L}(x_{k+1},\xi_{k+1})$. We assume that $\xi_{k+1}^i$, $i=1,...,m$ are i.i.d for fixed $k+1$. 
Also let the assumptions of the Theorem \ref{Th04} hold. Then
\begin{align}
& \E_{x_0, X_0, \dots, X_k} \varphi(y_k) - \varphi^* \leq \notag \\
& \leq \frac{C_1 LR_0^2V^2  }{k^{p}} + \frac{C_2 \sigma R_0 V }{\sqrt{mk}} +  C_3 k^{p-1}\delta, \notag
\end{align}
where $C_1=4 \sqrt{2}$, $C_2 = 16 \sqrt{2}$, $C_3 = 48$ and expectation is taken with respect to all randomness.
\label{Lm:amplMean}
\end{Lm}

\myproof 

1. Note that $d \left( \frac{x-x_0}{R_0} \right)$ is strongly convex with respect to the norm $\frac{1}{R_0}\|\cdot\|$ with parameter 1 and that the dual for this norm is the norm $R_0\|\cdot\|_*$.  Also note that with respect to the norm $\frac{1}{R_0}\|\cdot\|$ $(f_{\delta,L}(x),g_{\delta,L}(x))$ is the $(\delta, LR_0^2)$-oracle for $f(x)$. Also we have 
$$
\E_{\xi_{k+1}^1,\dots,\xi_{k+1}^m} \tilde{G}_{\delta,L}(x_{k+1}) = g_{\delta,L}(x_{k+1}),
$$
and
$$
\E_{\xi_{k+1}^1,\dots,\xi_{k+1}^m} R_0^2 \|\tilde{G}_{\delta,L}(x_{k+1})-g_{\delta,L}(x_{k+1})\|^2_* = \E_{\xi_{k+1}^1,\dots,\xi_{k+1}^m} R_0^2 \left\|\frac{1}{m}\sum_{i=1}^m{G_{\delta,L}(x_{k+1},\xi_{k+1}^i)}-g_{\delta,L}(x_{k+1})\right\|^2_* \stackrel{\reff{eq:varG}}{\leq}  \frac{\sigma^2 R_0^2 }{m}.
$$
Applying theorems \ref{Th02} and \ref{Th04} with changing $L$ to $LR_0^2$, $\sigma$ to $\frac{\sigma R_0 }{\sqrt{m}}$, $R$ to $V$ we obtain
\begin{align}
& \E_{x_0, X_0, \dots, X_k} \varphi(y_k) - \varphi^* \leq \notag \\
& \frac{\beta_k \E_{x_0} d \left( \frac{x^*-x_0}{R_0} \right) }{A_k} + \frac{\sum_{i=0}^k B_i \delta}{A_k} + \frac{1}{A_k}
 \sum_{i=0}^k \frac{B_i}{\beta_i-L}  \sigma^2 \leq \notag \\
& \leq \frac{C_1 LR_0^2 V^2}{k^{p}} + \frac{C_2 \sigma R_0 V }{\sqrt{mk}} +  C_3 k^{p-1}\delta. \notag 
\end{align}
\qed

Now we are ready to formulate the new algorithm for strongly convex problems.
\begin{algorithm}[h!]
\SetAlgoNoLine
\KwIn{The function $d(x)$, point $u_0$, number $R_0$ such that $\|u_0-x^*\| \leq R_0$, number $p \in [1,2]$.}
\KwOut{The point $u_{k+1}$.}

Set $k$ = 0.

Define 
\begin{equation}
N_k = \left\lceil \left(\frac{4eC_1LV^2}{\mu}\right)^{\frac{1}{p}}\right\rceil  .
\label{eq:SIGMA1N_0}
\end{equation}

\Repeat{}{
Define
\begin{align}
& m_k = \max \left\{1, \left\lceil \frac{16e^{k+2}  C_2^2 \sigma^2 V^2}{\mu^2 R_0^2 N_k } \right\rceil \right\},  \label{eq:SIGMA1m_k} \\
& R_k^2 = R_0^2 e^{-k} + \frac{2^p e C_3 \delta} {\mu(e-1)} \left(\frac{4eC_1LV^2}{\mu}\right)^{\frac{p-1}{p}} \left(1-e^{-k} \right). \label{eq:SIGMA1R_k}
\end{align}
Run Algorithm \ref{alg:SIGM} with $x_0=u_k$, prox-function $d \left( \frac{x-u_k}{R_k} \right)$ for $N_k$ steps using oracle $\tilde{G}^k_{\delta,L}(x) = \frac{1}{m_k}\sum_{i=1}^{m_k}{G_{\delta,L}(x,\xi^i)}$, where $\xi^i$, $i=1,...,m_k$ are i.i.d, on each step and sequences $\{\alpha_i\}_{i \geq 0}$, $\{\beta_i\}_{i \geq 0}$, $\{B_i\}_{i \geq 0}$ defined in Theorem \ref{Th04}. 

Set $u_{k+1}=y_{N_k}$, $k=k+1$.
    }
\caption{Stochastic Intermediate Gradient Method Accelerated (SIGMA)}
\label{alg:SIGMA1}
\end{algorithm}

\newpage
Let us prove the following result about rate of convergence of this algorithm.

\begin{Th}
After $k \geq 1$ outer iterations of the Algorithm \ref{alg:SIGMA1} we have
\begin{align}
&\E \varphi(u_{k}) - \varphi^* \leq \frac{\mu R_0^2}{2} e^{-k} + \frac{C_3 e 2^{p-1}}{e-1} \left(\frac{4eC_1LV^2}{\mu}\right)^{\frac{p-1}{p}} \delta, \label{eq:SIGMArate_phi} \\
& \E \|u_{k} - x^*\|^2 \leq R_0^2 e^{-k} + \frac{C_3 e 2^{p}}{\mu(e-1)} \left(\frac{4eC_1LV^2}{\mu}\right)^{\frac{p-1}{p}}  \delta. \label{eq:SIGMArate_x}
\end{align}

As a consequence if we choose error of the oracle $\delta$ satisfying
\begin{equation}
\delta \leq \frac{\e (e-1)}{2^pC_3e} \left(\frac{4eC_1LV^2}{\mu}\right)^{\frac{1-p}{p}}
\label{eq:SIGMADelta}
\end{equation}

then we need $N= \left\lceil \ln \left( \frac{\mu R_0^2}{ \e }\right) \right\rceil$ outer iterations and no more than
$$
\left(1+ \left(\frac{4eC_1LV^2}{\mu}\right)^{\frac{1}{p}} \right) \left( 1+ \ln \left( \frac{\mu R_0^2}{ \e }\right) \right) + \frac{16e^3C_2^2\sigma^2V^2}{\mu \e (e-1)}
$$
oracle calls to provide $\E \varphi(u_{N}) - \varphi^* \leq \e$. 
\label{Th06}
\end{Th}

\myproof Obviously \reff{eq:SIGMArate_x} follows from \reff{eq:SIGMArate_phi} and \reff{eq:PhiStrConv}. 
Let us prove the inequality
\begin{equation}
\E \varphi(u_{k}) - \varphi^* \leq \frac{\mu R_0^2}{2} e^{-k} + \frac{C_3 e 2^{p-1} }{e-1} \left(\frac{4eC_1LV^2}{\mu}\right)^{\frac{p-1}{p}}\left(1-e^{-k}\right) \delta
\label{eq:SIGMArate_phi1}
\end{equation}
for all $k \geq 1$. Obviously then we will have \reff{eq:SIGMArate_phi} as a consequence.
Let us prove \reff{eq:SIGMArate_phi1} for $k=1$. It follows from the Lemma \ref{Lm:amplMean} that
$$
\E \varphi(y_{N_0}) - \varphi^* \leq \frac{C_1 LR_0^2 V^2}{N_0^{p}} + \frac{C_2 \sigma R_0 V }{\sqrt{m_0N_0}} +  C_3 N_0^{p-1}\delta.
$$
From \reff{eq:SIGMA1N_0} we have
$$
\frac{C_1 LR_0^2 V^2}{N_0^{p}} \leq \frac{C_1 LR_0^2 V^2}{\frac{4eC_1LV^2}{\mu}} \leq \frac{\mu R_0^2}{4e}, \quad
C_3 N_0^{p-1} \delta \leq  \frac{C_3 e 2^{p-1}}{e-1} \left(\frac{4eC_1LV^2}{\mu}\right)^{\frac{p-1}{p}} \left(1-e^{-1}\right)   \delta.
$$
From \reff{eq:SIGMA1m_k} we have 
$$
\frac{C_2 \sigma R_0 V }{\sqrt{m_0N_0}} \leq \frac{C_2 \sigma R_0 V }{\sqrt{\frac{16e^{2}  C_2^2 \sigma^2 V^2}{\mu^2 R_0^2 N_0 } N_0}} \leq \frac{\mu R_0^2}{4e}.
$$
And we obtain \reff{eq:SIGMArate_phi1} for $k=1$.
Let us now assume that \reff{eq:SIGMArate_phi1} holds for $k=j$ and prove that it holds for $k=j+1$.
It follows from \reff{eq:SIGMArate_phi1} for $k=j$ and \reff{eq:PhiStrConv} that
$$
\E \|u_j-x^*\|^2\leq \frac{2}{\mu} \left(\E \varphi(u_j) - \varphi^* \right)\leq \frac{2}{\mu}  \left( \frac{\mu R_0^2}{2} e^{-j} + \frac{C_3 e 2^{p-1} }{e-1} \left(\frac{4eC_1LV^2}{\mu}\right)^{\frac{p-1}{p}}\left(1-e^{-j}\right) \delta \right) = R_j^2.
$$ 
After $N_{j}$ iterations of the Algorithm \ref{alg:SIGM} starting from the point $u_{j} = y_{N_{j-1}}$ applying Lemma \ref{Lm:amplMean} we have
$$
\E \varphi(y_{N_{j}}) - \varphi^* \leq \frac{C_1 LR_j^2 V^2}{N_{j}^{p}} + \frac{C_2 \sigma R_j V }{\sqrt{m_{j}N_{j}}} +  C_3 N_{j}^{p-1}\delta.
$$ 
From \reff{eq:SIGMA1N_0} we have
$$
\frac{C_1 LR_j^2 V^2}{N_{j}^{p}} \leq \frac{C_1 LR_{j}^2 V^2}{\frac{4eC_1LV^2}{\mu}} \leq \frac{\mu R_j^2}{4e}, \quad
C_3 N_{j}^{p-1} \delta \leq C_3 2^{p-1} \left(\frac{4eC_1LV^2}{\mu}\right)^{\frac{p-1}{p}}  \delta.
$$
From \reff{eq:SIGMA1m_k} we have 
$$
m_{j} \geq \frac{16e^{j+2}  C_2^2 \sigma^2 V^2}{\mu^2 R_0^2 N_{j} } \geq \frac{16e^{2}  C_2^2 \sigma^2 V^2}{\mu^2 N_{j} \left( R_0^2 e^{-j} + \frac{2^p e C_3 \delta} {\mu(e-1)} \left(\frac{4eC_1LV^2}{\mu}\right)^{\frac{p-1}{p}} \left(1-e^{-j} \right) \delta \right) } =  \frac{16e^{2}  C_2^2 \sigma^2 V^2}{\mu^2 R_{j}^2 N_{j} }, 
$$
and
$$
\frac{C_2 \sigma R_j V }{\sqrt{m_jN_j}} \leq \frac{C_2 \sigma R_j V }{\sqrt{\frac{16e^{2}  C_2^2 \sigma^2 V^2}{\mu^2 R_j^2 N_j } N_j}} \leq \frac{\mu R_j^2}{4e}.
$$
Finally we have
\begin{align}
&\E \varphi(u_{j+1}) - \varphi^* = \E \varphi(y_{N_{j}}) - \varphi^* \leq \frac{\mu R_j^2}{2e} + C_3 2^{p-1} \left(\frac{4eC_1LV^2}{\mu}\right)^{\frac{p-1}{p}} \delta = \notag \\
& = \frac{1}{e} \left( \frac{\mu R_0^2}{2} e^{-j} + \frac{C_3 e 2^{p-1} }{e-1} \left(\frac{4eC_1LV^2}{\mu}\right)^{\frac{p-1}{p}}\left(1-e^{-j}\right) \delta \right) + C_3 2^{p-1} \left(\frac{4eC_1LV^2}{\mu}\right)^{\frac{p-1}{p}}   \delta = \notag \\
&= \frac{\mu R_0^2}{2} e^{-(j+1)} + \frac{C_3 e 2^{p-1} }{e-1} \left(\frac{4eC_1LV^2}{\mu}\right)^{\frac{p-1}{p}}\left(1-e^{-(j+1)}\right) \delta. \notag
\end{align}
So we have obtained that \reff{eq:SIGMArate_phi1} holds for $k=j+1$ and by induction it holds for all $k \geq 1$.

If we choose $\delta$ satisfying \reff{eq:SIGMADelta} and perform $N=\left\lceil \ln \left( \frac{\mu R_0^2}{ \e }\right)\right\rceil$ outer iterations of SIGMA method we will obtain from \reff{eq:SIGMArate_phi} that
$$
\E \varphi(u_{N}) - \varphi^* \leq \frac{\mu R_0^2}{2} e^{-N} + \frac{C_3 e 2^{p-1}}{e-1} \left(\frac{4eC_1LV^2}{\mu}\right)^{\frac{p-1}{p}} \delta \leq \frac{\e}{2}+\frac{\e}{2} = \e.
$$
It remains to calculate the number of oracle calls to obtain an $\e$-solution $\E \varphi(u_{N}) - \varphi^* \leq \e$. We perform $N$ outer iterations (counting from $0$ to $N-1$) on each outer iteration $k$ we perform $N_k$ inner iterations and on each inner iteration we call the oracle $m_k$ times. So the total number of oracle calls is
\begin{align}
&{\mathcal C} (\e) = \sum_{k=0}^{N-1}{N_k m_k} \leq \sum_{i=0}^{N-1}{N_k \left(1 + \frac{16 C_2^2 e^2 \sigma^2 V^2 e^k}{\mu^2 R_0^2 N_k } \right)} \leq N N_0+{\frac{16 C_2^2 \sigma^2 V^2 e^{N}}{\mu^2 R_0^2 (e-1)}} \leq \notag \\
& \left(1+ \left(\frac{4eC_1LV^2}{\mu}\right)^{\frac{1}{p}} \right) \left( 1+ \ln \left( \frac{\mu R_0^2}{ \e }\right) \right) + \frac{16e^3C_2^2\sigma^2V^2}{\mu \e (e-1)} \leq \notag \\
& \leq \left(1+ \left(\frac{62LV^2}{\mu}\right)^{\frac{1}{p}} \right) \left( 1+ \ln \left( \frac{\mu R_0^2}{ \e }\right) \right) + \frac{96000\sigma^2V^2}{\mu \e }  . \notag
\end{align}
\qed

\subsection{Method with Bounded Large Deviations}
In this subsection we assume that the prox-function has quadratic growth with parameter $V^2$ with respect to the chosen norm:
\begin{equation}
d(x) \leq \frac{V^2}{2}\|x\|^2, \quad \forall x \in \R^n. 
\label{eq:dQuadrGrowth}
\end{equation}
Several examples of such prox-functions can be found in \cite{JudNest_2014}.

Now we present the modification of the SIGMA algorithm with a bound for large deviations.
\begin{algorithm}[h!]
\SetAlgoNoLine
\KwIn{The function $d(x)$, point $u_0$, number $R_0$ such that $\|u_0-x^*\| \leq R_0$, number $p \in [1,2]$, number $N \geq 1$ of outer iterations, confidence level $\Lambda$.}
\KwOut{The point $u_{N}$.}

Set $k$ = 0.

Define 
\begin{equation}
N_k = \left\lceil \left(\frac{6eC_1LV^2}{\mu}\right)^{\frac{1}{p}}\right\rceil  .
\label{eq:SIGMA2N_0}
\end{equation}

\Repeat{$k=N-1$}{
Define
\begin{align}
& m_k = \max \left\{1, \left\lceil \frac{36e^{k+2}  C_2^2 \sigma^2 V^2 \left(1+\ln \left(\frac{3N}{\Lambda} \right) \right)^2}{\mu^2 R_0^2 N_k } \right\rceil, \left\lceil \frac{144e^{k+2}  C_4^2 \sigma^2 \ln \left(\frac{3N}{\Lambda} \right) }{\mu^2 R_0^2 N_k } \right\rceil \right\},  \label{eq:SIGMA2m_k} \\
& R_k^2 = R_0^2 e^{-k} + \frac{2^p e C_3 \delta} {\mu(e-1)} \left(\frac{6eC_1LV^2}{\mu}\right)^{\frac{p-1}{p}} \left(1-e^{-k} \right), \label{eq:SIGMA2R_k} \\
& Q_k = \left\{ x \in Q : \|x-u_k\|^2 \leq R_k^2 \right\}. \label{eq:SIGMA2Q_k}
\end{align}
Run Algorithm \ref{alg:SIGM} applied to problem $\min_{x \in Q_k} \varphi(x)$ with $x_0=u_k$, prox-function $d \left( \frac{x-u_k}{R_k} \right)$ for $N_k$ steps using oracle $\tilde{G}^k_{\delta,L}(x) = \frac{1}{m_k}\sum_{i=1}^{m_k}{G_{\delta,L}(x,\xi^i)}$, where $\xi^i$, $i=1,...,m_k$ are i.i.d, on each step and sequences $\{\alpha_i\}_{i \geq 0}$, $\{\beta_i\}_{i \geq 0}$, $\{B_i\}_{i \geq 0}$ defined in Theorem \ref{Th04}. 

Set $u_{k+1}=y_{N_k}$, $k=k+1$.
    }
\caption{Stochastic Intermediate Gradient Method Accelerated 2}
\label{alg:SIGMA2}
\end{algorithm}

\newpage
Let us prove the following result about the rate of convergence of this algorithm.

\begin{Th}
After $N$ outer iterations of the Algorithm \ref{alg:SIGMA2} we have
\begin{align}
&\Prob \left\{ \varphi(u_{N}) - \varphi^* >  \frac{\mu R_0^2}{ 2 } e^{-N} + \frac{2^{p-1} e C_3 \delta} {(e-1)} \left(\frac{6eC_1LV^2}{\mu}\right)^{\frac{p-1}{p}} \delta \right\} \leq \Lambda \label{eq:SIGMA2eL}
\end{align}

As a consequence if we choose error of the oracle $\delta$ satisfying
\begin{equation}
\delta \leq \frac{\e (e-1)}{2^{p}C_3e} \left(\frac{6eC_1LV^2}{\mu}\right)^{\frac{1-p}{p}}
\label{eq:SIGMA2Delta}
\end{equation}
and choose $N= \left\lceil \ln \left( \frac{\mu R_0^2}{ \e }\right) \right\rceil$ outer iterations, then we will need no more than
\begin{align}
&\left(1+ \left(\frac{6eC_1LV^2}{\mu}\right)^{\frac{1}{p}} \right) \left( 1+ \ln \left( \frac{\mu R_0^2}{ \e }\right) \right) + \frac{36e^{3}  C_2^2 \sigma^2 V^2 }{\mu(e-1) \e } \left(1+\ln \left(\frac{3}{\Lambda}\left(1+\ln \left( \frac{\mu R_0^2}{ \e }\right) \right) \right) \right)^2 + \notag \\ 
& +\frac{144e^{3}  C_4^2 \sigma^2 }{\mu \e (e-1) } \ln \left(\frac{3}{\Lambda} \left(1+\ln \left( \frac{\mu R_0^2}{ \e }\right)\right) \right)  \label{eq:SIGMA2Compl}
\end{align}
oracle calls to provide $\Prob \{ \varphi(u_{N}) - \varphi^* > \e \} \leq \Lambda$. 
\label{Th07}
\end{Th}
\myproof

Let $A_k$, $k\geq0$ be event $A_k = \left\{\varphi(u_k) - \varphi^* \leq \frac{\mu R_k^2}{ 2 } \right\}$ and $\bar{A}_k$ be its complement. Let us prove first that for $k \geq 1$ 
\begin{equation}
\Prob \left\{\left. \varphi(u_{k}) - \varphi^* > \frac{\mu R_k^2}{ 2 } \right| A_{k-1} \right\} \leq \frac{\Lambda}{N}.
\label{eq:PrDeltaCond}
\end{equation}
Since the event $A_{k-1}$ holds we have from \reff{eq:PhiStrConv} that 
$$
\left\|u_{k-1} - x^*\right\|^2 \leq \frac{2}{\mu} \left(\varphi(u_{k-1}) - \varphi^* \right) \leq R_{k-1}^2.
$$
Hence the solution of the problem $\min_{x \in Q_{k-1}} \varphi(x)$ is the same as the solution of the initial problem \reff{eq:PrStateInit}. Let us denote $D_{k-1} = \max_{x,y\in Q_{k-1}} \|x-y\|$. Clearly $D_{k-1} \leq 2R_{k-1}$. Note that 
$ D_{k-1} = R_{k-1} \max_{x,y\in Q_{k-1}}\frac{\|x~-~y\|}{R_{k-1}}$ and the diameter of the set $Q_{k-1}$ with respect to the norm $\frac{\|\cdot\|}{R_{k-1}}$ is not greater than 2. Using the same argument as in the proof of the Lemma \ref{Lm:amplMean} but now using \reff{eq:dQuadrGrowth} and applying theorems \ref{Th03} and \ref{Th05} with changing $L$ to $LR_{k-1}^2$, $\sigma$ to $\frac{\sigma R_{k-1} }{\sqrt{m_{k-1}}}$, $R$ to $V$, $D$ to 2 we obtain
\begin{align}
& \Prob \left\{ \varphi(u_k) - \varphi^* > \left. \frac{C_1 LR_{k-1}^2V^2 }{N_{k-1}^{p}} + \frac{C_2( 1+ \Omega )\sigma R_{k-1} V }{\sqrt{m_{k-1}N_{k-1}}} +  C_3 N_{k-1}^{p-1}  \delta + \frac{2C_4 R_{k-1} \sigma \sqrt{\Omega} }{\sqrt{m_{k-1}N_{k-1}}} \right| A_{k-1}  \right\} \leq \frac{\Lambda}{N}  , \label{eq:Th06pr01}
\end{align}
where $C_1=4 \sqrt{2}$, $C_2 = 16 \sqrt{2}$, $C_3 = 48$, $C_4 = 4 \sqrt{3}$, $\Omega = \ln \left(\frac{3N}{\Lambda} \right)$.

From \reff{eq:SIGMA2N_0} we have
$$
\frac{C_1 LR_{k-1}^2 V^2}{N_{k-1}^{p}} \leq \frac{C_1 LR_{k-1}^2 V^2}{\frac{6eC_1LV^2}{\mu}} \leq \frac{\mu R_{k-1}^2}{6e}, \quad
C_3 N_{k-1}^{p-1} \delta \leq C_3 2^{p-1} \left(\frac{6eC_1LV^2}{\mu}\right)^{\frac{p-1}{p}}  \delta.
$$
From \reff{eq:SIGMA2m_k} we have 
\begin{align}
&m_{k-1} \geq \frac{36e^{k+1}  C_2^2 \sigma^2 V^2 (1+\Omega)^2}{\mu^2 R_0^2 N_{k-1} } \geq \frac{36e^{2}  C_2^2 \sigma^2 V^2  (1+\Omega)^2 }{\mu^2 N_{k-1} \left( R_0^2 e^{-(k-1)} + \frac{2^p e C_3 \delta} {\mu(e-1)} \left(\frac{6eC_1LV^2}{\mu}\right)^{\frac{p-1}{p}} \left(1-e^{-(k-1)} \right) \delta \right) } = \notag \\
& =  \frac{36e^{2}  C_2^2 \sigma^2 V^2  (1+\Omega)^2}{\mu^2 R_{k-1}^2 N_{k-1} }, \notag
\end{align}
and
$$
\frac{C_2 \sigma R_{k-1} V  (1+\Omega)}{\sqrt{m_{k-1}N_{k-1}}} \leq \frac{C_2 \sigma R_{k-1} V  (1+\Omega)}{\sqrt{\frac{36e^{2}  C_2^2 \sigma^2 V^2  (1+\Omega)^2}{\mu^2 R_{k-1}^2 N_{k-1} } N_{k-1}}} \leq \frac{\mu R_{k-1}^2}{6e}.
$$
Also from \reff{eq:SIGMA2m_k} we have 
$$
m_{k-1} \geq \frac{144e^{k+1}  C_4^2 \sigma^2 \Omega}{\mu^2 R_0^2 N_{k-1} } \geq \frac{144e^{2}  C_4^2 \sigma^2 \Omega}{\mu^2 N_{k-1} \left( R_0^2 e^{-(k-1)} + \frac{2^p e C_3 \delta} {\mu(e-1)} \left(\frac{6eC_1LV^2}{\mu}\right)^{\frac{p-1}{p}} \left(1-e^{-(k-1)} \right) \delta \right) } =  \frac{144e^{2}  C_4^2 \sigma^2 \Omega}{\mu^2 R_{k-1}^2 N_{k-1} }, 
$$
and
$$
\frac{2C_4 R_{k-1} \sigma \sqrt{\Omega} }{\sqrt{m_{k-1}N_{k-1}}} \leq \frac{2C_4 R_{k-1} \sigma \sqrt{\Omega} }{\sqrt{\frac{144e^{2}  C_4^2 \sigma^2 \Omega}{\mu^2 R_{k-1}^2 N_{k-1} } N_{k-1}}} \leq \frac{\mu R_{k-1}^2}{6e}.
$$
Finally we have
\begin{align}
& \frac{C_1 LR_{k-1}^2V^2 }{N_{k-1}^{p}} + \frac{C_2( 1+ \Omega )\sigma R_{k-1} V }{\sqrt{m_{k-1}N_{k-1}}} +  C_3 N_{k-1}^{p-1}  \delta + \frac{2C_4 R_{k-1} \sigma \sqrt{\Omega} }{\sqrt{m_{k-1}N_{k-1}}} \leq \notag \\
& \leq \frac{\mu R_{k-1}^2}{2e} + C_3 2^{p-1} \left(\frac{6eC_1LV^2}{\mu}\right)^{\frac{p-1}{p}} \delta = \notag \\
& = \frac{1}{e} \left( \frac{\mu R_0^2}{2} e^{-(k-1)} + \frac{C_3 e 2^{p-1} }{e-1} \left(\frac{6eC_1LV^2}{\mu}\right)^{\frac{p-1}{p}}\left(1-e^{-(k-1)}\right) \delta \right) + C_3 2^{p-1} \left(\frac{6eC_1LV^2}{\mu}\right)^{\frac{p-1}{p}}   \delta = \notag \\
&= \frac{\mu R_0^2}{2} e^{-k} + \frac{C_3 e 2^{p-1} }{e-1} \left(\frac{6eC_1LV^2}{\mu}\right)^{\frac{p-1}{p}}\left(1-e^{-k}\right) \delta = \frac{\mu R_k^2}{2}. \notag
\end{align}
Hence from \reff{eq:Th06pr01} we have \reff{eq:PrDeltaCond}. 

Also for all $k = 1, \dots, N$ we have
\begin{align}
& \Prob\left\{ \varphi(u_{k}) - \varphi^* > \frac{\mu R_k^2}{ 2 } \right\} = \Prob \left\{\left. \varphi(u_{k}) - \varphi^* > \frac{\mu R_k^2}{ 2 }\right| A_{k-1} \cup \bar{A}_{k-1} \right\} = \notag \\
& = \Prob \left\{\left. \varphi(u_{k}) - \varphi^* > \frac{\mu R_k^2}{ 2 } \right| A_{k-1}  \right\} \Prob \{ A_{k-1} \} + \Prob \left\{\left. \varphi(u_{k}) - \varphi^* > \frac{\mu R_k^2}{ 2 } \right| \bar{A}_{k-1}  \right\} \Prob \{\bar{A}_{k-1} \} \stackrel{\reff{eq:PrDeltaCond}}{\leq} \notag \\
& \leq \frac{\Lambda}{N} + \Prob \{\bar{A}_{k-1} \} =  \frac{\Lambda}{N} + \Prob \left\{\varphi(u_{k-1}) - \varphi^* > \frac{\mu R_{k-1}^2}{ 2 }\right\}.\notag
\end{align}
Using that $\Prob\{A_0\}=1$ and summing up these inequalities we obtain
$$
\Prob\left\{ \varphi(u_{N}) - \varphi^* > \frac{\mu R_0^2}{ 2 } e^{-N} + \frac{2^{p-1} e C_3 \delta} {(e-1)} \left(\frac{6eC_1LV^2}{\mu}\right)^{\frac{p-1}{p}} \delta \right\} \leq \Prob\left\{ \varphi(u_{N}) - \varphi^* > \frac{\mu R_N^2}{ 2 } \right\} \leq \Lambda.
$$
Making the same arguments as in the proof of the Theorem \ref{Th06} we obtain the complexity bound \reff{eq:SIGMA2Compl}.
\qed

\section{Conclusion and discussion}
We have proposed SIGM which can be used for convex composite optimization problems with stochastic inexact oracle. This method has rate of convergence $\Theta\left(\frac{LR^2}{k^p} + \frac{\sigma R}{\sqrt{k}} + k^{p-1}\delta\right)$. Also we have provided bounds for large deviation for the error of the method $\varphi(y_k) - \varphi^*$ which has the same asymptotic dependence on $k$. This method also provides several degrees of freedom for adapting it to the problem at hand.
\begin{enumerate}
	\item Depending on the relations between error of the oracle $\delta$ and constant $L$ we can choose the value of $p \in [1,2]$ to have optimal trade-off between error accumulation and rate of convergence.
	\item We can introduce randomization to the problem if stochastic approximation of the gradient is cheaper to obtain than the real gradient. Since the rate of convergence depends only on $k$ - number of iterations but not on the number of calls of the oracle, we can use Monte Carlo idea and generate several realizations of stochastic approximation of the gradient on each iteration. This can reduce the variation of the stochastic approximation from $\sigma^2$ to $\sigma^2/m$, where $m$ is the number of generated realizations of $G(x,\xi)$.
	\item The notion of $(\delta,L)$-oracle allows to use the proposed method to solve non-smooth problems. It was shown in \cite{DGN_FOM_2011} that convex non-smooth function with H\"older continuous subgradient $\|g(x) -g(y) \|_* \leq L_{\nu} \|x-y\|^{\nu}, \quad \nu \in [0,1]$ can be equipped for any $\delta >0$ with $(\delta,L)$-oracle, where $L=L_{\nu} \left[ \frac{L_{\nu}(1-\nu)}{2 \delta (1+\nu)} \right]^{\frac{1-\nu}{1+\nu}}$.
	\item Since the method uses general prox-function and norm we can choose them optimally depending on the geometry of the problem. For example, if $Q$ is a standard simplex in $n$-dimensional space, $d(x) = - \ln n + \sum_{i=1}^n x_i \ln x_i$, $h(x) =0$, then the optimization on the steps like \reff{eq:algzk}, \reff{eq:alghxkp1} can be done explicitly \cite{NemYud_1983}.
	\item The method allows to solve composite optimization problems such as LASSO $\|A x -b \|_2^2 + \lambda \|x\|_1 \to \min$.
	\item If we know that the function $\varphi(x)$ is strongly convex, we can use restart technique to have better rate of convergence. In this case we also have a modification of the method which provide an $(\e,\Lambda)$-solution $\hat{u}$ satisfying $\Prob \{ \varphi(\hat{u}) - \varphi^* > \e \} \leq \Lambda$.
\end{enumerate}

\section{Acknowledgments}
Authors would like to thank professor Yurii Nesterov and professor Arkadi Nemirovski for useful discussions.



\addcontentsline{toc}{section}{References}

\begin{thebibliography}{9}

\bibitem{Evt_1982}
Yu. Evtushenko. \emph{Methods of Solving Extremal Problems and Their Application in Optimization Systems}.
Moscow: Nauka, 1982.

\bibitem{Polyak_1987}
B.T. Polyak. \emph{Introduction to Optimization}. Optimization Software Inc, 1987

\bibitem{NemYud_1983} A. Nemirovski, D. Yudin. \emph{Problem complexity and method efficiency in optimization}. Wiley Interscience Series in Discrete Mathematics. John Wiley, XV, 1983.

\bibitem{Nest_2004}
Y. E. Nesterov. \emph{Introductory Lectures on Convex Optimization: a basic course}. Kluwer Academic
Publishers, Massachusetts, 2004.

\bibitem{KhaTarErl_1988}
L.Khachiyan, S.Tarasov, and E.Erlich. \emph{ The inscribed ellipsoid method}.
Soviet Math. Dokl. (In Russian) , 298 (1988).

\bibitem{NemNes_1994}
A. Nemirovski and Yu.Nesterov. \emph{Interior point polynomial methods in
convex programming: Theory and Applications}. SIAM, Philadelphia, 1994.

\bibitem{Nes_2012}
Yu. Nesterov. \emph{Subgradient methods for huge-scale optimization problems}.
CORE Discussion Paper, 2, 2012

\bibitem{Tib_1996}
Tibshirani, R. \emph{Regression shrinkage and selection via the lasso}. Journal of the Royal Statistical Society, Series B 58 (1): 267Ц288, 1996.

\bibitem{DGN_IGM_2013} O. Devolder, F. Glineur and Yu. Nesterov. \emph{Intermediate Gradient Methods for Smooth Convex Problems with Inexact Oracle}. CORE Discussion Paper 2013/17, available at
\url{http://www.uclouvain.be/cps/ucl/doc/core/documents/coredp2013_17web.pdf}.

\bibitem{DGN_FOM_2011}
O. Devolder, F. Glineur and Yu. Nesterov. \emph{First-order Methods of Smooth
Convex Optimization with Inexact Oracle}. CORE Discussion Paper 2011/2, available at \url{http://www.optimization-online.org/DB_FILE/2010/12/2865.pdf}.

\bibitem{JudLanNemShap_2009}
A. Juditsky, G. Lan, A. Nemirovski, A. Shapiro \emph{Stochastic approximation approach to stochastic programming}. SIAM Journal on Optimization. 2009, 19(4), pp. 1574--1609.

\bibitem{LanNemShap_2012}
G. Lan, A. Nemirovski and A. Shapiro. \emph{Validation analysis of mirror descent stochastic approximation method.} Mathematical Programming Serie A, 2012, 134(2), pp. 425--458.

\bibitem{Dev_PhD_2013}
O. Devolder. \emph{Exactness, Inexactness and Stochasticity in First-Order Methods for Large-Scale Convex Optimization}, PhD thesis (2013).

\bibitem{GhLan_SCO1_2012}
S.Ghadimi, G.Lan, \emph{Optimal Stochastic Approximation Algorithms for Strongly Convex Stochastic Composite Optimization I: A Generic Algorithmic Framework}, SIAM J. Optim., 2012, 22(4), pp. 1469--1492.

\bibitem{GhLan_SCO2_2013}
S.Ghadimi, G.Lan, \emph{Optimal Stochastic Approximation Algorithms for Strongly Convex Stochastic Composite Optimization II: Shrinking Procedures and Optimal Algorithms}, SIAM J. Optim., 2013, 23(4), pp. 2061--2089.

\bibitem{JudNest_2014}
A. Juditsky, Yu. Nesterov, \emph{Primal-dual subgradient methods for minimizing uniformly convex
functions.}, 2014, available at \url{http://arxiv.org/abs/1401.1792}.

%
%
%

\end{thebibliography}
\end{document}